\input amstex
\input amsppt.sty
\documentstyle{amsppt}

\hsize=5.3in 
\vsize=7.5in 
\magnification=\magstep1 
\nologo

\input epsf         
\epsfverbosetrue    

\def\E{{\Bbb E}}
\def\V{\text{Vol}}

\def\p{\bold p = (\bold p_1, \dots, \bold p_N)}
\def\q{\bold q = (\bold q_1, \dots, \bold q_N)}
\def\Picture #1 by #2 (#3){ 
  \vbox to #2{
  \hrule width #1 height 0pt depth 0pt \vfill \special{picture #3}}}
\def\scaledpicture #1 by #2 (#3 scaled #4){{
 \dimen0=#1 \dimen1=#2
 \divide\dimen0 by 1000 \multiply\dimen0 by #4 
 \divide\dimen1 by 1000 \multiply\dimen1 by #4 
 \Picture \dimen0 by \dimen1 (#3 scaled #4)}} 

\comment
\vskip12pt
\centerline{\scaledpicture 5.03in by 2.65in (Fig6.1.1 scaled 694)}
\centerline{\bf Figure 6.1.1} 
\endcomment

\topmatter 
\title Pushing disks apart - The Kneser-Poulsen conjecture in the plane 
\endtitle 
\author  K\'aroly Bezdek and Robert Connelly \endauthor
\affil E\"otv\"os University and Cornell University \endaffil
\address
E\"otv\"os University 
Department of Geometry 
H-1053 Budapest 
Kecskem\'eti utca 10-12
Hungary
\endaddress
\email kbezdek\@ludens.elte.hu\endemail
\address
Department of Mathematics
Malott Hall
Cornell University 
Ithaca, NY 14853
USA
\endaddress
\email connelly\@math.cornell.edu \endemail
\thanks The first author was partially supported by the Hung. Nat. 
Sci. Found. (OTKA), grant no. T029786. \endthanks
\date July 31, 2001 \enddate
\abstract
 We give a proof of the planar case of a longstanding conjecture of
Kneser (1955) and Poulsen (1954). In fact, we prove more by showing
that if a finite set of disks in the plane is rearranged so that the
distance between each pair of centers does not decrease, then the area
of the union does not decrease, and the area of the intersection does
not increase.
\endabstract 
\endtopmatter

\document

\subhead 1. Introduction \endsubhead
 
If $\p$ and $\q$ are two configurations of $N$ points, where each
$\bold p_i \in \E^n$ and each $\bold q_i \in \E^n$ is such that for
all $1 \le i < j \le N$, $|\bold p_i -\bold p_j| \le |\bold q_i -
\bold q_j|$, we say that $\bold q$ is an {\it expansion} of $\bold p$
(and $\bold p$ is a {\it contraction} of $\bold q$).  If $\bold q$ 
is an expansion of $\bold p$, then there
may or may not be a continuous motion $\bold p(t) = (\bold p_1(t),
\dots, \bold p_N(t))$ for $0 \le t \le 1$ such that $\bold p(0) =\bold
p$ and $\bold p(1) = \bold q$, and $|\bold p_i(t) - \bold p_j(t)|$ is
monotone increasing.  When there is such a motion, we say that $\bold
q$ is a {\it continuous expansion} of $\bold p$.  Let $B(\bold p_i,
r_i)$ be the closed ball of radius $r_i \ge 0$ in $\E^n$ about the
point $\bold p_i$, $\V_n$ represent the $n$-dimensional volume, and
$|\dots|$ be the Euclidean norm, so $|\bold p_i - \bold p_j|$ is the
Euclidean distance between $\bold p_i$ and $\bold p_j$.

In 1954 Poulsen [21] and in 1955 Kneser [18] independently 
conjectured the following for the case when $r_1 = \dots = r_N$:

\proclaim{Conjecture \rom{1}} If $\q$ is an expansion
of $\p$ in $\E^n$, then
$$ 
\V_n \left[\bigcup_{i=1}^N B(\bold p_i, r_i)\right] \le 
\V_n \left[\bigcup_{i=1}^N B(\bold q_i, r_i)\right]. \tag1
$$
\endproclaim

 We will prove this conjecture for the case of the 
plane, $n=2$, and with the same hypothesis the following related conjecture,
which was mentioned in [17] by Klee and Wagon.

\proclaim{Conjecture \rom{2}} If $\q$ is an expansion
of $\p$ in $\E^n$, then
$$ 
\V_n \left[\bigcap_{i=1}^N B(\bold p_i, r_i)\right] \ge 
\V_n \left[\bigcap_{i=1}^N B(\bold q_i, r_i)\right]. \tag2
$$
\endproclaim

In [5] Bollob\'as proved Conjecture 1, for $n=2$, when $r_1 = \dots =
r_N$ and $\bold q$ is a continuous expansion of $\bold p$.  In [4]
Bern and Sahai proved Conjecture 1 and Conjecture 2 for $n=2$, 
but again with the additional assumption that $\bold
q$ is a continuous expansion of $\bold p$.  Then in [11] Csik\'os
extended Bollob\'as's result to arbitrary radii for $n=2$, and later
in [10] Csik\'os proved Conjecture (1) under the assumption that $\bold q$ 
is a continuous expansion of $\bold p$. 
 In [7] Capoyleas showed (2) for congruent 
radii in the plane, but assuming that $\bold q$ is a continuous
expansion of $\bold p$.  In [13] Gromov proved (2) for arbitrary
radii, but only for $N \le n+1$. Then in [8] Capoyleas and Pach
proved (1) for arbitrary radii in all dimensions, but again only for $N
\le n +1$.  In these cases, it is not hard to show that if $\bold q$
is an expansion of $\bold p$, then it is a continuous expansion, a
property that does not hold even for $n+2$ points in $\E^n$.  In all
of these cases, it is assumed or implicitly holds that the
configuration $\bold q$ is a continuous expansion of $\bold p$.

In the following, we will use a formula of Csik\'os describing the
derivative of the volume of the union of $4$-dimensional balls, when
their centers are expanding smoothly, to show that the area of the
union of $2$-dimensional disks increases when one configuration of
centers is an expansion of another, even when there is no continuous
expansion in the plane.  See Section 6 of this paper for a related
result which is that a particular weighted surface volume changes
monotonically under smooth continuous expansions.  For such smooth
continuous expansions this result extends the first result of Csik\'os
in [11].

Conjecture 1, with all the radii equal, was repeated by Hadwiger in
[14].  Later it was included in a list of problems by Valentine in
[25], Klee in [16], Croft, Falconer, and Guy in [9], Moser and
Pach in [20], and Klee and Wagon in [17], mentioning, in particular,
the case of disks in the plane.  This is the case that we prove here.

\subhead 2. Connecting configurations in higher dimensions \endsubhead

Our plan is to use results about continuous (or differentiable)
motions of configurations of points in a higher dimension to get
information about pairs of configurations in a lower dimension.  The
following lemma, which is fairly well-known, is essentially the same
as formula (8) in Alexander [1], where the $\sqrt{1-t}$ and $\sqrt{t}$
in Alexander's formula is replaced by $\cos(\pi t)$ and $\sin(\pi t)$,
respectively, and this is composed with a rotation to bring the final
image back to the original copy of $\E^n$.  See Gromov [13], and
Capolyeas and Pach [8], for a related result with a different
proof. This allows us to connect configurations in a higher
dimension.  We regard $\E^n$ as the subset $\E^n = \E^n \times \{\bold
0\} \subset \E^n \times \E^n = \E^{2n}$.

\proclaim{Lemma \rom{1}} Suppose that $\p$ and 
$\q$ are two configurations in $\E^n$.  Then the following  
is a continuous motion $\bold p(t) = (\bold p_1(t), \dots, \bold p_N(t))$ 
in $\E^{2n}$, 
that is analytic in $t$, such that $\bold p(0)=\bold p$, 
$\bold p(1)=\bold q$ 
and for $0 \le t \le 1$,
$|\bold p_i(t)-\bold p_j(t)|$ is monotone:
$$
\bold p_i(t) = \left({\bold p_i + \bold q_i \over 2} + 
(\cos{\pi t}) {\bold p_i - \bold q_i \over 2}, 
(\sin{\pi t}) {\bold p_i - \bold q_i \over 2}\right), 
\quad 0 \le i < j \le N. 
\tag3
$$

\endproclaim

\demo{Proof} We calculate:
$$\align
4 \,|\bold p_i(t) - \bold p_j(t)|^2  =& \,
|(\bold p_i - \bold p_j) - (\bold q_i - \bold q_j)|^2 +
|(\bold p_i - \bold p_j) + (\bold q_i - \bold q_j)|^2 \\ &+ 
2(\cos{\pi t})(|\bold p_i - \bold p_j|^2 - |\bold q_i - \bold q_j|^2)
.
\endalign$$
This function is monotone, as required.
\enddemo

As stated here, the distances $|\bold p_i(t)-\bold p_j(t)|$ could be
monotone increasing or decreasing, but we will only need the case when
$\bold q$ is an expansion of $\bold p$ and thus all distances are
monotone increasing.  (Of course, we regard the constant function as
monotone.)

\subhead 3.  Main results \endsubhead

We say that a configuration $\q$ is a piecewise-smooth expansion of
$\p$ if $\bold q$ is a continuous expansion of $\bold p$, and all the
coordinates of all the points are infinitely differentiable functions
of the parameter $t$ except for a finite number of values of $t$.  The
following theorem and its corollaries are our main results. 
Also we regard $\E^n$ as the subset $\E^n = \E^n \times \{\bold
0\} \subset \E^n \times \E^2 = \E^{n+2}$.

\proclaim{Theorem \rom{1}} Let $\p$ and 
$\q$ be two configurations in $\E^n$ such that 
$\bold q$ is a piecewise-smooth expansion of $\bold p$ in $\E^{n+2}$.  
Then the conclusions (1) and (2)
of Conjecture 1 and Conjecture 2 hold in $\E^n$.
\endproclaim

The proof of this result will occupy the next few sections. 
The following includes the Kneser-Poulsen conjecture in the plane.

\proclaim{Corollary \rom{1}}  Let $\p$ and 
$\q$ be two configurations in $\E^2$ such that 
$\bold q$ is an arbitrary expansion of $\bold p$. 
Then (1) and (2) hold for $n = 2$.
\endproclaim

\demo{Proof} Apply Lemma 1 to the configurations 
$\bold p$ and $\bold q$ to
get that $\bold q$ is an analytic expansion of $\bold p$ in $\E^4$.  Then 
Theorem 1 applies, and the area inequalities follow.
\enddemo

The following is obtained by taking the limit as $r \rightarrow
\infty$ in Corollary 1, where $r_1 = \dots = r_N =r$. 
It is one of the main results in [24] of Sudakov, in [1] of
Alexander, and in [8] of Capoyleas and Pach. Although all these papers
do not prove Corollary 1, it is
explained carefully in [8] how to derive Corollary 2 from the 
Kneser-Poulsen conjecture in the plane.

\proclaim{Corollary \rom{2}} If $\q$ is an arbitrary expansion of $\p$ 
in $\E^2$, then the length of the perimeter of the convex hull of
$\bold p$ is less than or equal to the length of the perimeter of the
convex hull of $\bold q$.
\endproclaim

The following is an immediate consequence of Theorem 1 and formula
(3).

\proclaim{Corollary \rom{3}} If $\q$ is an arbitrary expansion of 
$\p$ in $\E^n$, and the vectors $\bold p_i - \bold q_i$, for all 
$1\le i \le N$, lie in a
$2$-dimensional subspace of $\E^n$, then both (1) and (2) hold.
\endproclaim

The following related version of Theorem 1 follows from its proof.
Right after the proof of Theorem 1 in Section 7 we will mention the
adjustments for the proof of Remark 1.
\vskip .05in  

\noindent {\bf Remark 1.} Let $\p$ and 
$\q$ be two configurations in $\E^n$ such that for some integer $m$,
$\bold q$ is a piecewise-smooth expansion of $\bold p$ in $\E^m$, where the
expansion is given by $\bold p(t)$, $\bold p(0) = \bold p$ and $\bold
p(1) = \bold q$, but the dimension of the affine span of $\bold p(t) =
(\bold p_1(t), \dots, \bold p_N(t))$ is at most $(n+2)$-dimensional
and is piecewise-constant.  Then the conclusions (1) and (2) of
Conjecture 1 and Conjecture 2 hold in $\E^n$.

\vskip .05in
The following generalizes a result of Gromov in [13], who proved it 
in the case $N \le n+1$.

\proclaim{Corollary \rom{4}} If $\q$ is an arbitrary expansion of 
$\p$ in $\E^n$, and $N \le n+3$, then then both (1) and (2) hold.
\endproclaim
\demo{Proof} Apply Lemma 1, to get the analytic expansion $\bold p(t)$
for $0 \le t \le 1$ between $\bold p$ and $\bold q$.  By taking the
determinant of an appropriate number of coordinates of an appropriate
subset of the vectors $\bold p_i(t) - \bold p_j(t)$ it follows that
the dimension of the affine span of $\bold p(t)$ is
piecewise-constant.  By assumption, the $n+3$ points can have an
affine span of dimension no larger than $n+2$.  Then Remark 1 applies.
\enddemo

As an example of how one might apply Theorem 1 in higher dimensions to 
expansions that are not continuous, we present the following result.
\proclaim{Corollary \rom{5}}Let $\q$ be an  expansion of $\p$ in $\E^n$ 
such that for some $\lambda > 1$, for each $i=1,\dots,N$ either 
$\bold q_i = \bold p_i$ or $\bold q_i = \lambda \bold p_i$.  Then 
$\bold q$ is a smooth expansion of $\bold p$ in $\E^{n+1}$ and thus 
(1) and (2) hold for any $r_i >0$, for $i=1,\dots,N$.
\endproclaim
\demo{Proof} In the definition of the motion, we replace (3) with 
$$
\bold p_i(t) = \left({\bold p_i + \bold q_i \over 2} + 
(\cos{\pi t}) {\bold p_i - \bold q_i \over 2}, 
(\sin{\pi t}) \left|{\bold p_i - \bold q_i \over 2}\right|\right), 
1 \le i \le N, 
$$
which lives in $\E^{n+1}$.  To check that it is an expansion, consider 
any $1 \le i < j \le N$.  If either $\bold q_i = \bold p_i$ or 
$\bold q_j = \bold p_j$, then a calculation similar to the one in 
the proof of Lemma 1 applies, and $|\bold p_i(t) - \bold p_j(t)|$ is
monotone increasing.  Otherwise, consider the case when 
$\bold q_i = \lambda \bold p_i$ and $\bold q_j = \lambda \bold p_j$. 
We calculate
$$
{d \over dt}|\bold p_i(t) - \bold p_j(t)|^2 =
{\pi \over 2}(\lambda - 1)^2 (\sin{\pi t}) |\bold p_i - \bold p_j|^2
\left[{\lambda +1 \over \lambda -1}  + (\cos{\pi t})
\left(\left({ |\bold p_i| - |\bold p_j|\over |\bold p_i - \bold p_j|} 
\right)^2 -1\right) \right].
$$
This is non-negative, which implies that $|\bold p_i(t) - \bold p_j(t)|$ 
is monotone increasing.
\enddemo

See Section 8 for an example of this sort of expansion.  
For a different approach to the case when 
the configurations are similar see [6] by Bouligand.  When 
the sets forming the union are not spherical balls, but translates of 
a convex  set, then in [22], Rehder shows that the volume of the union 
does not decrease when the sets are dilated.  However, a general 
expansive rearrangement of convex sets different from ellipsoids 
can have the volume of the intersection (union) 
increase (decrease), as shown in [19] by Meyer, Reisner and 
Schmuckenschl\"ager. 

\subhead 3.  Nearest point and farthest point Voronoi diagrams \endsubhead
  
For a given configuration $\p$ of points in $\E^n$, and radii 
$r_1, \dots, r_N$, consider the following sets:
$$\align
C_i &= \{\,\bold p_0 \in \E^n \mid \text{for all}\,\, 
j \ne i,\,\, |\bold p_0 - \bold p_i|^2 -r_i^2 \le 
|\bold p_0 - \bold p_j|^2 -r_j^2 \,\}\\
C^i &= \{\,\bold p_0 \in \E^n \mid \text{for all}\,\, 
j \ne i,\,\, |\bold p_0 - \bold p_i|^2 -r_i^2 \ge 
|\bold p_0 - \bold p_j|^2 -r_j^2 \,\}
\endalign$$

The set $C_i$ is the closed {\it extended nearest point Voronoi region} of 
points $\bold p_0$ whose {\it power}, 
$|\bold p_0~-~\bold p_i|^2  - r_i^2$ with respect to $\bold p_i$
is less than or equal to the power of any other point of the 
configuration with respect to $\bold p_i$.  
There is a good
discussion of how this decomposition fit into
the kind of problems that we are considering in [12] by Edelsbrunner.
The set $C^i$ is often called the {\it extended farthest point Voronoi
region} of points and there is a good discussion of this 
in [23] by Seidel.

We now restrict each of the sets by intersecting them with a ball of
radius $r$ centered at $\bold p_i$.
$$\align
C_i(r) &= C_i \cap B(\bold p_i, r)\\
C^i(r) &= C^i \cap B(\bold p_i, r)
\endalign$$

We will call the sets $\{C_i\}_{i=1}^N$ and $\{C^i\}_{i=1}^N$ the {\it
nearest and farthest point Voronoi cells}, respectively.  We shall be
interested especially in the collections $\{C_i(r_i)\}_{i=1}^N$ and
$\{C^i(r_i)\}_{i=1}^N$, which we call the {\it nearest and
farthest point truncated Voronoi cells}.  Figure 1 shows some examples
of these sets in the plane.

\centerline{
\epsfysize=2.0in                   
\epsfbox{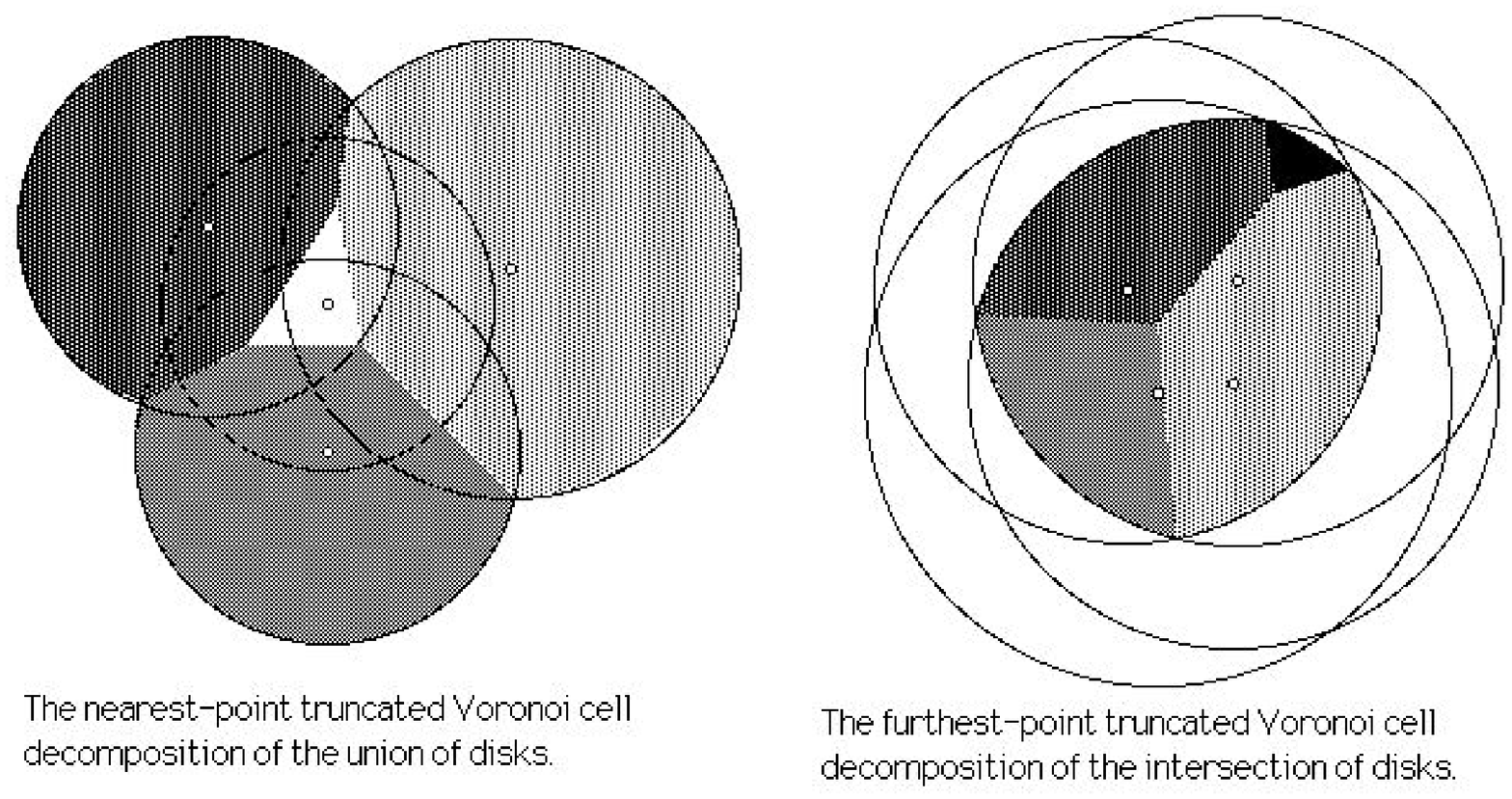}
}
\midspace{0.1in} \caption{Figure 1}

For each $i \ne j$ let $W_{ij} = C_i \cap C_j$ and $W^{ij} = C^i \cap
C^j$, and for any $\bold p_0 \in \E^n$ and $r > 0$, define 
$W_{ij}(\bold p_0, r) = W_{ij} \cap B(\bold p_0, r)$ 
and $W^{ij}(\bold p_0, r) = W^{ij} \cap B(\bold p_0, r)$.  
  
The $W_{ij}$ and $W^{ij}$ are
the {\it walls} between the nearest point and farthest point Voronoi
cells. Note that some of the walls may be empty or low-dimensional,
but in any case the walls lie in a hyperplane of dimension $n-1$.  
Define the following for $\bold r = (r_1, r_2, \dots, r_N)$:
$$
X_n(\bold p, \bold r) = \bigcup_{i=1}^N B(\bold p_i, r_i),\qquad
X^n(\bold p, \bold r) = \bigcap_{i=1}^N B(\bold p_i, r_i).
$$
We
need the following easily verified properties of these Voronoi diagrams.  
We will use 
$\text{Bdy}[X]$ to denote the boundary of a set $X$ in $\E^n$.
\roster 
\item"{\bf (i.)}"$\{C_i(r_i)\}_{i=1}^N$ is a tiling of 
$X_n(\bold p, \bold r)$ and $\{C^i(r_i)\}_{i=1}^N$ is a 
tiling of $X^n(\bold p, \bold r)$. 
\item"{\bf (ii.)}"$\text{Bdy}[X_n(\bold p, \bold r)] 
\cap B(\bold p_i, r_i) = 
\text{Bdy}[X_n(\bold p, \bold r)] \cap C_i(r_i)$ and \newline
$\text{Bdy}[X^n(\bold p, \bold r)] \cap B(\bold p_i, r_i) = 
\text{Bdy}[X^n(\bold p, \bold r)] \cap C^i(r_i)$.
\item"{\bf (iii.)}" $W_{ij}(\bold p_i, r_i)  =
W_{ij}(\bold p_j, r_j)$ and 
$W^{ij}(\bold p_i, r_i)  =
W^{ij}(\bold p_j, r_j) $.
\item"{\bf (iv.)}" When $W_{ij}(\bold p_i, r_i) \ne \emptyset$, the vector 
$\bold p_j - \bold p_i$ is a positive scalar multiple of 
the outward pointing normal to the boundary of $C_i(r_i)$ at 
$W_{ij}(\bold p_i, r_i)$.  Similarly, when 
$W^{ij}(\bold p_i, r_i) \ne \emptyset$, the vector 
$\bold p_j - \bold p_i$ is a negative scalar multiple of 
the outward pointing normal to the boundary of $C^i(r_i)$ at 
$W_{ij}(\bold p_i, r_i)$.
\endroster
So we get the following.

\proclaim{Lemma \rom{2}} For $r \le s$, $W_{ij}(\bold p_i, r) 
\subseteq W_{ij}(\bold p_i, s)$, 
and  $W^{ij}(\bold p_i, r) \subseteq W^{ij}(\bold p_i, s)$.
\endproclaim

\subhead 4. Integral formulas \endsubhead

One of the key ideas to prove Theorem 1 is a relation between the
surface volume of the union (and intersection) of the higher-dimensional
balls and the area of the union (and intersection) of lower dimensional
disks.  First we state a lemma from calculus.

\proclaim{Lemma \rom{3}}  Let $X$ be a compact integrable set in 
$\E^{n+2}$ that is a solid of revolution about $\E^n$.  In other 
words the projection of 
$X\cap\{\E^n \times \{(s\cos{\theta},s\sin{\theta})\}$ 
into $\E^n$ is an integrable set $X(s)$ independent of $\theta$. Then
$$
\V_{n+2} \left[X\right] = 2\pi \int_0^{\infty} \V_n \left[X(s)\right] s\, ds.
$$
\endproclaim

We specialize to the case when the set $X$ is the  
intersection of a ball of radius $r$, and half-spaces whose boundary is
orthogonal to $\E^n$. 

In the following $\bold p$ is a configuration of points in
$\E^n \subset \E^{n+2}$.  We are especially interested in the relation of 
the  volume 
of $C_i(r)=C_i(r,n)$ and $C^i(r)=C^i(r,n)$ in $\E^n$ to the volume of the 
corresponding truncated Voronoi cell $C_i(r,n+2)$ and $C^i(r,n+2)$ 
in $\E^{n+2}$.
\proclaim{Lemma \rom{4}} If $\bold p$ is a configuration of points in
$\E^n \subset \E^{n+2}$, then
$$\align
\V_{n+2} \left[C_i(r,n+2)\right] &= 
2\pi \int_0^r \V_n \left[C_i(s,n)\right]s\, ds\\
\V_{n+2} \left[C^i(r,n+2)\right] &= 
2\pi \int_0^r \V_n \left[C^i(s,n)\right]s\, ds
\endalign$$
\endproclaim
\demo{Proof} It is clear, in both cases, that  
$C_i(r,n+2)$ and $C^i(r,n+2)$ are 
compact sets of revolution. Let $B^{n+2}(\bold p_i, r)$ denote the closed 
ball of radius $r$ in $\E^{n+2}$.  Then 
$B^{n+2}(\bold p_i, r)\cap \{\E^n \times \{(s\cos{\theta},s\sin{\theta})\}$ 
is an $n$-dimensional ball 
of radius $\sqrt{r^2-s^2}$ in a translate of $\E^n$, and thus by Lemma 3 
 we have that 
$$
\V_{n+2} \left[C_i(r,n+2)\right] = 
2\pi \int_0^r \V_n \left[C_i(\sqrt{r^2-s^2},n)\right]s\, ds.
$$
But if we make the change of variable $u=\sqrt{r^2-s^2}$, we get the 
desired integral.  A similar calculation works for $C^i(r,n+2)$.
\enddemo
Applying the fundamental theorem of calculus we get the following.

\proclaim{Corollary \rom{6}} We have 
$$\align
{d \over dr}\V_{n+2} \left[C_i(r,n+2)\right] &= 
2\pi r \V_n \left[C_i(r,n)\right], \text{and} \\ 
{d \over dr}\V_{n+2} \left[C^i(r,n+2)\right] &= 
2\pi r \V_n \left[C^i(r,n)\right].
\endalign$$
\endproclaim
\vskip 0.05in
\noindent{\bf Remark 2.} We can interpret 
${d \over dr}\V_{n+2} \left[C_i(r,n+2)\right]$, evaluated at
$r=r_i$, as the 
$(n+1)$-dimensional surface volume of   
$\text{Bdy}[X_{n+2}(\bold p,\bold r)] \cap C_i(r_i,n+2)$,  since the
radius vector for the ball is orthogonal to that part of the boundary of
$C_i(r_i,n+2)$.  We can make a similar identification for 
$\text{Bdy}[X^{n+2}(\bold p, \bold r)] \cap C^i(r_i,n+2)$.

\subhead 5. Csik\'os's formula \endsubhead

Suppose that $\bold p(t) = (\bold p_1(t), \dots, \bold p_N(t))$, for 
$0 \le t \le 1$, is a 
smooth motion of the configuration $\bold p = \bold p(0)$ in some 
Euclidean space $\E^{n}$.  Let $d_{ij}= |\bold p_i(t) - \bold p_j(t)|$, 
and let $d'_{ij}$ be the $t$-derivative of $d_{ij}$.  Then  
 Csik\'os's Theorem 4.1,  
for unions of balls, in [10] is the following.  For intersections of balls,
we indicate the appropriate adjustments.

\proclaim{Theorem 2} Let $n \ge 2$ and let $\bold p(t)$ be a smooth 
motion of a configuration in $\E^{n}$ such that for each $t$, the 
points of the configuration are pairwise distinct.  Then regarding the 
following as functions of $t$,
$V_n(t,\bold r) = \text{Vol}_n[X_n(\bold p(t),  \bold r)]$ and 
$V^n(t,\bold r) = \text{Vol}_n[X^n(\bold p(t),  \bold r)]$ are 
differentiable and,
$$\align
{d \over dt} V_n(t,\bold r) &= \sum_{1 \le i < j \le N} 
d'_{ij} \text{Vol}_{n-1}[W_{ij}(\bold p_i(t), r_i)],\\
{d \over dt} V^n(t,\bold r) &= \sum_{1 \le i < j \le N} 
-d'_{ij} \text{Vol}_{n-1}[W^{ij}(\bold p_i(t), r_i)].
\endalign$$
\endproclaim
\demo{Proof} For the case of unions of balls, this is the same result as
in [10].  For the case of intersections, the proof proceeds in a very 
similar way, but when one uses property (iv.), there is a sign change.
\enddemo

The following is a result in [15] of Kirszbraun.  There are other
simple elementary proofs, for example in [17] as described by Klee and
Wagon and described by Alexander in [2].  It is immediate from
Theorem 2 and Lemma 1, which was also pointed out by Alexander.

\proclaim{Corollary \rom{7}} If the configuration $\bold p$ is a contraction 
of the configuration $\bold q$ in $\E^n$, and 
$\bigcap_{i=1}^N B(\bold q_i, r_i)$ is non-empty, then 
$\bigcap_{i=1}^N B(\bold p_i, r_i)$ is non-empty as 
well.
\endproclaim

\subhead 6. Expanding the configuration \endsubhead

In order to get a global relation between the $(n+1)$-dimensional volume 
of the surface of our sets in $\E^{n+2}$ and the $n$-dimensional volume 
of our sets in $\E^n$, we consider a particular 
deformation of just the radii, fixing the configuration $\bold p$.  For 
each $i = 1, 2, \dots, N$ and $0 \le s$, define $r_i(s)=\sqrt{r_i^2 + s}$.  
Each $r_i$ is constant, and the function $r_i(s)$ is independent of the 
parameter $t$. We assume that each $r_i > 0$.  Then we calculate that
$$
{d \over ds}r_i(s) = {1 \over 2 r_i(s)}. \tag4
$$
Now define $\bold r(s) = (r_1(s), \dots, r_N(s))$, and regard 
$\V_n[X_n(\bold p(t),\bold r(s))]=V_n(t,s)$ and 
$\V_n[X^n(\bold p(t),\bold r(s))]=V^n(t,s)$ as 
functions of both variables $t$ and $s$.  Throughout we assume that all
$r_i > 0$.

\proclaim{Lemma \rom{5}} Let the motion of the configuration $\bold p(t)$
be at least continuously differentiable.  Then the functions 
$V_n(t,s)$ and $V^n(t,s)$ are
continuously differentiable in $t$ and $s$ simultaneously, and for 
fixed $t$, the extended nearest point and farthest point Voronoi 
cells are constant.
\endproclaim
\demo{Proof} Recall that a point $\bold p_0$ is in an extended Voronoi cell
$C_i$ or $C^i$, when for all $j \ne i$,  
$|\bold p_0 - \bold p_i|^2 - |\bold p_0 - \bold p_j|^2 - r_i(s)^2 + r_j(s)^2$ 
is non-positive for $C_i$ and non-negative for $C^i$. But 
$r_i(s)^2 - r_j(s)^2 = r_i^2 - r_j^2$ is constant.  So each 
$C_i$ and  $C^i$ is constant as a function of $s$.  

Since $\bold p(t)$ is continuously differentiable, then the partial 
derivatives of $V_n(t,s)$ and $V^n(t,s)$ with respect to $t$ exist and 
are continuous.  Since the $(n-1)$-dimensional surface volume of 
boundaries of $X_n(\bold p, \bold r(s))$
and  $X^n(\bold p, \bold r(s))$ are continuous functions of $s$, the partial 
derivatives of $V_n(t,s)$ and $V^n(t,s)$ with respect to $s$ exist and
are continuous.  Thus $V_n(t,s)$ and $V^n(t,s)$ are both continuously 
differentiable with respect to $t$ and $s$ simultaneously. 
\enddemo

Note that we now can interchange the order of partial differentiation 
with respect to the variables $t$ and $s$.  Combining Lemma 5 and Theorem 2,
we get the following.

\proclaim{Lemma \rom{6}} Let $\bold p(t)$ be a smooth motion of a  
configuration in $\E^n$ such that for each $t$ the points of the 
configuration are pairwise distinct. Then the following hold:
$$\align
{\partial^2 \over \partial t \partial s} V_n(t,s) &= \sum_{1 \le i < j \le N} 
d'_{ij}{\partial \over \partial s}
\text{Vol}_{n-1}[W_{ij}(\bold p_i(t), r_i(s))],\\
{\partial^2 \over \partial t \partial s} V^n(t,s) &= \sum_{1 \le i < j \le N} 
-d'_{ij} {\partial \over \partial s} 
\text{Vol}_{n-1}[W^{ij}(\bold p_i(t), r_i(s))].
\endalign$$
Hence if $\bold p(t)$ is expanding, then by Lemma 2, 
${\partial  \over \partial s}V_n(t,s)$ is monotone 
increasing in $t$, and ${\partial  \over \partial s}V^n(t,s)$ 
is monotone 
decreasing in $t$.
\endproclaim

Bear in mind that we can replace 
$W_{ij}(\bold p_i(t), r_i(s))$ by $W_{ij}(\bold p_j(t), r_j(s))$ in the
 terms above by property (iii.).

Let $K_i(\bold p, \bold r)$ and $K^i(\bold p, \bold r)$ be the 
$(n-1)$-dimensional surface volume of 
$\text{Bdy}[X_n(\bold p, \bold r)]\cap C_i$ and
$\text{Bdy}[X^n(\bold p, \bold r)]\cap C^i$ respectively.  Then we 
observe the following, using Property (ii.).

\proclaim{Theorem 3} We can interpret ${\partial  \over \partial s}V_n(t,s)$ 
and  ${\partial  \over \partial s}V^n(t,s)$ evaluated at 
$\bold r = \bold r(0)$
as 
$${1 \over 2}\sum_{i=1}^N K_i(\bold p, \bold r)/r_i \quad \text{and} 
\quad {1 \over 2}\sum_{i=1}^N K^i(\bold p, \bold r)/r_i,
$$
the weighted $(n-1)$-dimensional volume of the boundary of 
$X_n(\bold p,\bold r)$ 
and $X^n(\bold p,\bold r)$
respectively. Thus under smooth expanding motions, these boundary 
volumes are monotone functions.
\endproclaim

For smooth motions this generalizes the result in [5] of Bollob\'as
for the plane as well as Csik\'os's other proof in [11] for not
necessarily congruent disks in the plane.  See Section 8 for comments,
however.

\subhead 7.  Proof of Theorem 1 \endsubhead

We now specialize to the case when the configuration is in $\E^n$, but the 
motion occurs in $\E^{n+2}$.  So we wish to connect the volumes of 
$\text{Vol}_n[X_n(\bold p,  \bold r(s))]= V_n(\bold p, \bold r(s))$ and 
$\text{Vol}_n[X^n(\bold p,  \bold r(s))]= V^n(\bold p, \bold r(s))$ in 
$\E^n$ to the corresponding volumes 
$V_{n+2}(\bold p, \bold r(s))$ and $V^{n+2}(\bold p, \bold r(s))$ in 
$\E^{n+2}$.
\proclaim{Lemma \rom{7}} Let $\p$ be a fixed configuration in 
$\E^n \subset \E^{n+2}$.  Then 
$$\align
{d \over ds}V_{n+2}(\bold p, \bold r(s)) &= \pi V_n(\bold p, \bold r(s)), 
\quad \text{and} \\
{d \over ds}V^{n+2}(\bold p, \bold r(s)) &= \pi V^n(\bold p, \bold r(s)).
\endalign$$
\endproclaim

\demo{Proof} By Property (i.), $V_{n+2}(\bold p, \bold r(s)) = \sum_{i=1}^N 
\V_{n+2} \left[C_i(r_i(s),n+2)\right]$; applying Corollary 6, the 
chain rule, and (4) we have that 
$$\align
{d \over ds}V_{n+2}(\bold p, \bold r(s)) &=  
\sum_{i=1}^N  {d \over ds}V_{n+2} (C_i(r_i(s),n+2))\\
&=\sum_{i=1}^N  {d \over dr_i(s)}V_{n+2} (C_i(r_i(s),n+2))
{dr_i(s) \over ds}\\
&=\sum_{i=1}^N  2 \pi r_i(s) V_{n} (C_i(r_i(s),n))
\left({1 \over 2r_i(s)}\right)\\
&=\pi V_n(\bold p, \bold r(s)).
\endalign$$
Similarly ${d \over ds}V^{n+2}(\bold p, \bold r(s))$ is calculated.
\enddemo

We are now in a position to show our main result.

\demo{Proof of Theorem 1}Suppose that the configuration $\q$ is 
an expansion of the configuration $\p$ in $\E^n$.  
By assumption, there 
is a piecewise-smooth expansion  
$\bold p(t) = (\bold p_1(t), \dots, \bold p_N(t))$, for $0 \le t \le 1$,
in $\E^{n+2}$ such that $\bold p(0) = \bold p$, and $\bold p(1) = 
\bold q$. So there is a finite number of sub-intervals of $[0,1]$, 
where each pair of 
points is distinct or remain coincident.  In either case, in the interior 
of each interval, 
by Lemma 6 applied to $\E^{n+2}$, we conclude that 
${d \over ds}V_{n+2}(\bold p(t), \bold r(s))$ is increasing in $t$.  By taking 
limits as $t$ approaches the endpoints of each interval, we have that 
${d \over ds}V_{n+2}(\bold p(t), \bold r(s))$ is increasing for all 
$0 \le t \le 1$.
Applying  Lemma 7,  
$\pi V_n(\bold p(0),\bold r(s)) = 
{d \over ds}V_{n+2}(\bold p(0), \bold r(s)) 
\le 
{d \over ds} V_{n+2}(\bold p(1), \bold r(s)) =
\pi V_n(\bold p(1),\bold r(s))$. 
Evaluating when $s=0$, we get the desired result.  A similar argument 
shows that $V^n(\bold p(0),\bold r) \ge V^n(\bold p(1),\bold r)$.
\enddemo  

\demo{Proof of Remark 1} Here the motion of the configuration 
$\bold p(t)$ is in $\E^m$, but the dimension of the affine span is at
most $n+2$ and the the dimension of the span of $\bold p(t)$ is
piecewise-constant. On each interval, while the dimension is constant,
it is possible to continuously, smoothly define an orthonormal
coordinate system, whose dimension is the dimension of the affine span
of the configuration $\bold p(t)$.  If the dimension of the affine
span is less than $n+2$, define additional coordinates so that there
is always an $(n+2)$-dimensional coordinate system during the interior
of each of the time intervals. For sufficiently small subintervals of
these intervals, the proof of Theorem 1 applies to these coordinate
systems.  So the $(n+1)$-dimensional weighted volume of the boundary
changes monotonically as before.  Then Lemma 7 applies, and we get the
desired result.
\enddemo

\subhead 8.  Examples and comments \endsubhead
Theorem 3 is delicate. If the
configuration $\bold q$ is an expansion of $\bold p$ but not a
continuous expansion, then even in the plane with disks of the same
radius, the length of the boundary of the union of disks  
may not be larger for $\bold q$ than for $\bold p$.  The example 
in Figure 2, due to Habicht and Kneser, 
described in [17],  shows this in the plane. 

\centerline{
\epsfysize=2.0in                   
\epsfbox{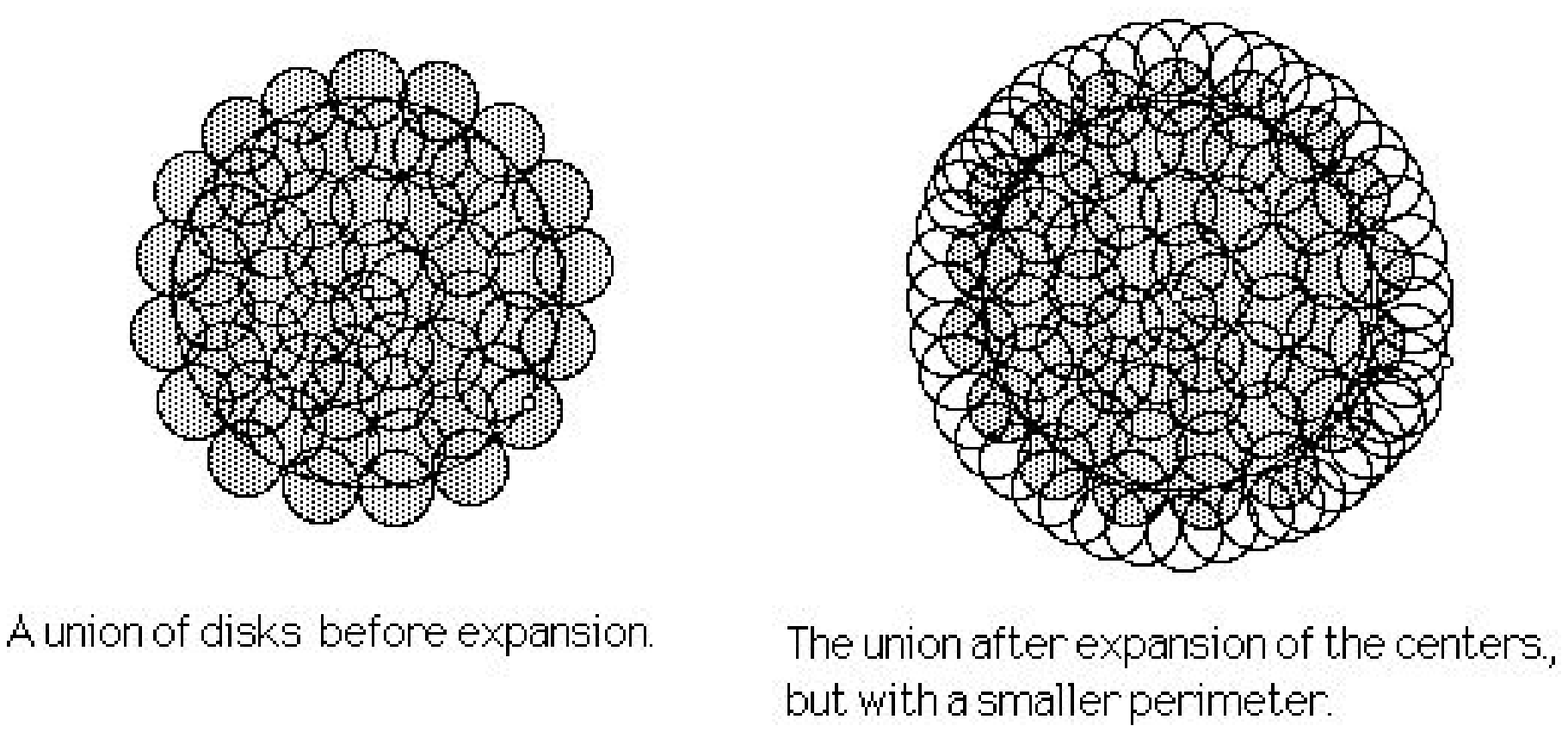}
}
\midspace{.1in} \caption{Figure 2}
Except for a small portion of its boundary, the inner shaded region 
is covered by a large number of congruent disks.  Then for a large $k$, 
there are $k$ disks that are arranged on the boundary as indicated.  This 
is the configuration $\bold p$.  Then some of the inner disks are moved 
radially outward covering almost all of the old boundary, leaving behind 
enough disks to still almost cover the original union.  This is the 
expanded configuration $\bold q$, and the associated disks almost cover the 
boundary of the disks about $\bold p$, but now the boundary is almost 
a perfect circle.  The ratio of the length of the boundary of the union of 
the disks about $\bold q$ to the length of the boundary of the union of 
the disks about $\bold p$ approaches $\pi/2 > 1$.  We do not know how to 
get a better ratio in the plane.  This example extends to higher 
dimensions.

If we have incongruent disks in the plane and a smooth motion, it can 
happen that the (unweighted) length of the boundary of the union can 
decrease while
the configuration is expanding.  The following example is very similar to 
the one described in [3] by Bern.

\centerline{
\epsfysize=2.0in                   
\epsfbox{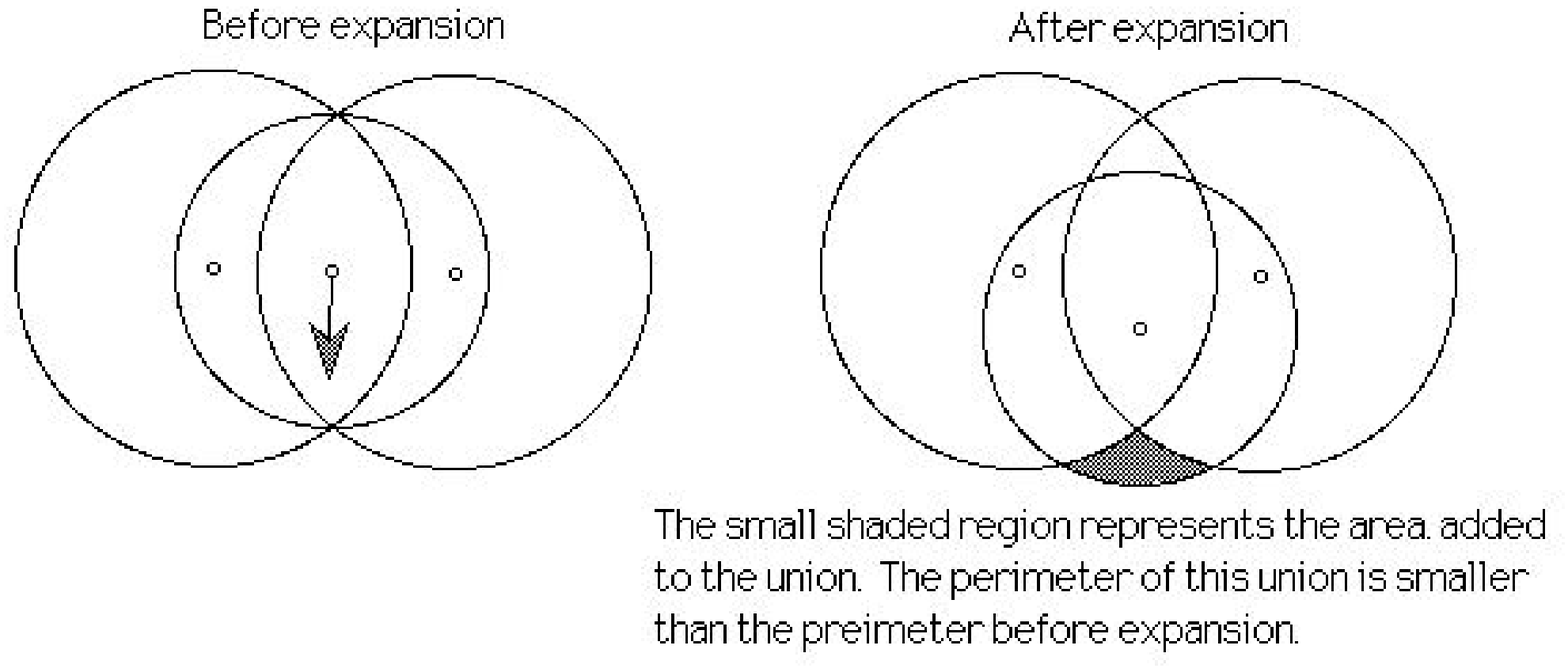}
}
\midspace{.1in} \caption{Figure 3}
The smaller disk moves as indicated, which is clearly a smooth
expansion of the three centers.  The shaded triangular region $\Delta$,
whose shape is close to an actual triangle, represents the additional
area. Two of the sides of $\Delta$ vanish as part of the boundary of
the union of the disks, and the third side is the new part of the
boundary.  The triangle inequality implies that the length of the
boundary decreases as the smaller circle moves.  Note that two of the
circles are the same length, and by choosing the two equal circles
sufficiently close to each other, this example will work when the
radius of third disk is arbitrarily close to the radius of the other
two.

We are still left with the question as to what happens for volumes of
expansions of unions and intersections of balls in higher dimensions.
It is possible that an extension of Lemma 6 could help.  We finish
with the following comment.
\vskip 0.05 in
\noindent{\bf Remark 3.} If the following inequalities hold for $k$ for a 
sufficiently smoothly expanding configuration  
$\bold p(t)$ in $\E^{4k}$, then 
Conjecture 1 and Conjecture 2 hold in $\E^{2k}$.
$$\align
{\partial^{k+1} \over \partial t (\partial s)^k} V_n(t,s) 
&= \sum_{1 \le i < j \le N} 
d'_{ij}{(\partial)^k \over (\partial s)^k}
\text{Vol}_{n-1}[W_{ij}(\bold p_i(t), r_i(s))] \ge 0,\\
{\partial^{k+1} \over \partial t (\partial s)^k}  V^n(t,s) 
&= \sum_{1 \le i < j \le N} 
-d'_{ij} {(\partial)^k \over (\partial s)^k} 
\text{Vol}_{n-1}[W^{ij}(\bold p_i(t), r_i(s))] \le 0.
\endalign$$

\enddocument